\documentclass[VANCOUVER,STIX1COL]{WileyNJD-v2}
\usepackage{moreverb}

\newcommand\BibTeX{{\rmfamily B\kern-.05em \textsc{i\kern-.025em b}\kern-.08em
		T\kern-.1667em\lower.7ex\hbox{E}\kern-.125emX}}

\articletype{Research Article}%

\received{<day> <Month>, <year>}
\revised{<day> <Month>, <year>}
\accepted{<day> <Month>, <year>}

\usepackage{amsmath,amsfonts,amssymb,mathtools}
\usepackage{url}
\usepackage{latexsym}
\usepackage{graphicx,overpic}

\newcommand{\R}{\mathbb R}
\newcommand{\bA}{\mathbf A}
\newcommand{\bB}{\mathbf B}
\newcommand{\bC}{\mathbf C}
\newcommand{\bL}{\mathbf L}

\newcommand{\bI}{\mathbf I}

\newcommand{\bP}{\mathbf P}

\newcommand{\bM}{\mathbf M}
\newcommand{\bQ}{\mathbf Q}
\newcommand{\bS}{\mathbf S}

\newcommand{\bV}{\mathbf V}

\newcommand{\bn}{\mathbf n}
\newcommand{\be}{\mathbf e}
\newcommand{\bp}{\mathbf p}

\newcommand{\bu}{\mathbf u}
\newcommand{\blf}{\mathbf f}
\newcommand{\bU}{\mathbf U}
\newcommand{\bv}{\mathbf v}

\newcommand{\bw}{\mathbf w}

\newcommand{\bx}{\mathbf x}

\newcommand{\T}{\mathcal T}

\newcommand{\tr}{{\rm tr}}

\newcommand{\divG}{{\mathop{\,\rm div}}_{\Gamma}}

\newcommand{\gradG}{\nabla_{\Gamma}}

\newcommand{\OGamma}{\Omega^\Gamma_h}

\newcommand{\nablaG}{\nabla_{\Gamma}}

\newcommand{\vect}[1]{\boldsymbol{\mathbf{#1}}}
\newcommand{\la}{\left\langle}
\newcommand{\ra}{\right\rangle}

\usepackage{subcaption}
\usepackage{graphicx}
\graphicspath{{img/}} 
\newcommand{\includegraphicsw}[2][1.]{\includegraphics[width=#1\linewidth]{#2}}
\usepackage{multirow}
\usepackage{hhline}
\usepackage{float} 

\DeclareMathOperator{\Curl}{curl}
\DeclareMathOperator{\vCurl}{\vect{curl}}

\newcommand{\LInfSpace}[1][\Gamma]{{L^\infty\left({#1}\right)}}

\newcommand{\rev}[1]{{\color{black}#1}}

\begin{document}

\title{Recycling augmented Lagrangian preconditioner in an incompressible fluid solver\protect\thanks{This
			work was  partially supported by US National Science Foundation (NSF) through grants DMS-1953535 and DMS-2011444.}}

\author[1]{Maxim A. Olshanskii*}

\author[2]{Alexander Zhiliakov}

\authormark{M.A. Olshanskii and A. Zhiliakov}

\address[1]{	\orgdiv{Department of Mathematics}, \orgname{University of Houston}, \orgaddress{ \state{Houston, TX}, \country{USA} } }

\address[2]{	\orgdiv{Department of Mathematics}, \orgname{University of Houston}, \orgaddress{ \state{Houston, TX}, \country{USA} } }

\corres{*Maxim A. Olshanskii, Department of Mathematics, University of Houston, Houston,TX 77204, USA. \email{maolshanskiy@uh.edu}}

\abstract[Abstract]{The paper discusses a reuse of matrix factorization as a building block in the Augmented Lagrangian (AL) and modified AL preconditioners for  non-symmetric saddle point  linear algebraic systems. The strategy is applied to solve two-dimensional incompressible fluid problems with efficiency rates independent of the Reynolds number. The solver is then tested to simulate  motion of a surface fluid, an example of a 2D flow motivated by an interest in lateral fluidity of inextensible viscous membranes.
	Numerical examples include the Kelvin--Helmholtz instability problem posed on the sphere and on the torus. \rev{Some new eigenvalue estimates for the AL preconditioner are derived.}}


\jnlcitation{\cname{%
		\author{Olshanskii M.A.} and
		\author{A. Zhiliakov} } (\cyear{2021}), 
	\ctitle{Recycling augmented Lagrangian preconditioner in an incompressible fluid solver}, \cjournal{Numer. Linear Algebra
		Appl.}, \cvol{...}.}

\maketitle

\section{Introduction}
Augmented Lagrangian (AL) preconditioning is a potent technique  that has been developed to solve some highly non-symmetric algebraic systems having a saddle point structure~\cite{benzi2006augmented,benzi2011modified,borm2012,de2007two,farrell2020augmented,he2011augmented,he2018combining,heister2013efficient,moulin2019augmented,olshanskii2008augmented,ur2008comparison}. The need to treat such problems numerically emerges from the discretization of systems of PDEs describing the motion of incompressible viscous fluid with dominating inertia effects.
Adopting the terminology of fluid mechanics, the AL approach augments the velocity subproblem of the system using a suitably weighted incompressibility constraint. This leads to a well conditioned pressure Schur complement matrix, but makes the velocity submatrix more difficult to solve or to precondition. Already in the original work~\cite{benzi2006augmented} a special multigrid method has been used to overcome the difficulty associated with preconditioning the velocity block, and recently this technique  was extended in~\cite{farrell2020reynolds,farrell2019augmented}. Nevertheless, the specialized multilevel approach is efficient only if a hierarchy of nested discretizations is available and only for certain finite element velocity--pressure pairs.
In the present paper we advocate a more general but still efficient, way to handle the velocity subproblem in the AL approach.  The proposed method consists of computing a (possibly incomplete) LU factorization of the velocity block (or velocity sub-blocks in a modified AL approach) with further recycling of the factors over several  time steps. The factorization can be updated when the velocity field variations significantly change the transport term  in the equations. A simpler strategy adopted here consists of updating preconditioner when the number of FGMRES iterations exceeds a threshold.  We shall see that for realistic unsteady 2D flows this results in a very efficient approach, which is robust with respect to the Reynolds number and calls for only a small number of full factorizations over a long-time simulation.

Employing matrix factorizations in algebraic solvers for equations governing the flow of viscous incompressible fluids is not a new theme. It is standard to factorize the discrete pressure Poisson equation. More recently, studies were done regarding different strategies to perform incomplete LU factorization of the coupled systems for velocity and pressure~\cite{dahl1992ilu,konshin2017lu,konshin2015ilu,segal2010preconditioners}. We note that the latter cannot be done with the help of position-based ILU, since the pressure block in the matrix may be zero. The augmented Lagrangian approach provides a framework to apply factorization only to the velocity matrix, while retaining the overall excellent preconditioning properties. The velocity matrix results from the discretization of  an elliptic part of the system. Therefore it is typically a positive definite matrix and LU factorization is stable without any preprocessing.

Large scale 3D simulations lead to algebraic systems which are still too expensive to factorize exactly and alternative ways of treating the velocity submatrix, e.g. based on geometric/algebraic multigrid, domain decomposition methods or incomplete factorization, can be more feasible and practical. The situation is different for 2D problems, where acceptable resolution is often achieved using the number of  degrees of freedom  affordable by  state-of-the-art direct solvers executed on a desktop machine. Traditionally 2D flows have been considered as a mathematical idealization of real-life 3D phenomena. However, recently we see a growing interest in understanding and solving fluid systems posed on 2D surfaces (see, e.g.,  \cite{bonito2020divergence,brandner2020finite,fries2018higher,gross2019meshfree,Jankuhn1,jankuhn2019higher,lederer2020divergence,nitschke2019hydrodynamic,nitschke2012finite,
olshanskii2018finite,olshanskiiinf,olshanskii2019penalty,reuther2018solving})
 as physically motivated continuum based  models of thin material layers exhibiting lateral viscosity such as lipid bilayers and plasma membranes~\cite{GurtinMurdoch75,rangamani2013interaction,torres2019modelling}.

This recent interest in models describing lateral fluidity of material surfaces motivates our choice of the test fluid problem here. We consider the  Navier--Stokes equations posed on a geometrically steady surface. These equations govern the tangential viscous motions of the surface fluid  subject to the inextensibility condition.
A geometrically unfitted discretization technique known as the trace finite element method (Trace FEM)~\cite{ORG09,olshanskii2016trace} is applied to handle the systems numerically. The augmentation is added on finite element (FE) level in the form  of the grad-div stabilization~\cite{olshanskii2002low,olshanskii2004grad}.
\rev{For such a setting, we prove eigenvalue bounds for the preconditioned system that extend a result from \cite{benzi2006augmented}
for the case of the unfitted surface FEM and the  FE-level augmentation. Compared to the case treated in \cite{benzi2006augmented} and other publications, the FE-level augmentation delivers a different algebraic structure. The latter required us to search for different arguments to prove the eigenvalue bounds.}
A preconditioned iterative method with recycled factorizations is then applied to solve the linearized Navier--Stokes equations on each time step of numerical simulations. The whole approach is used to compute two interesting surface flows:   the Kelvin--Helmholtz instability problem on a sphere and a torus. We observe a notable difference in the evolution of (large) vortices in these two settings.

\rev{Summarizing, the paper contributes to the field of developing fast and reliable solvers for fluid problems by (i) Introducing a simple and efficient strategy to reuse factorizations of positive definite matrices leading to a solver robust with respect to the Reynolds number;
(ii) Extending the AL preconditioner and its eigenvalue analysis to the trace FEM discretization of the linearized surface Navier--Stokes equations. Our construction can be useful for other unfitted FEMs for fluid equations.
Equipped with these tools we simulate for the first time the Kelvin--Helmholtz instability on a torus.
}

The remainder of the paper consists of three sections. In section~\ref{sectioncont} we recall the surface Navier--Stokes equations and apply $\vect P_2$--$P_1$ Trace FEM to discretize them. Section~\ref{s:matr} considers the properties of the resulting linear algebraic systems, introduces a preconditioner, \rev{eigenvalue bounds},  and recycling strategy. A proof of the bounds is given in the Appendix section. Finally, section~\ref{s_num} collects and discusses the results of numerical experiments.

\section{Surface Navier--Stokes problem and its discretization} \label{sectioncont}
Consider a clo\-sed sufficiently smooth surface  $\Gamma$ embedded in $\mathbb{R}^3$.  For $\bn$, the outward pointing unit normal on $\Gamma$, the orthogonal projection on the tangential plane is given by $\bP=\bP(\bx)\coloneqq \bI - \bn(\bx)\bn(\bx)^T$, $\bx \in \Gamma$.  Following~\cite{Jankuhn1} we formulate the surface fluid equation in terms of tangential differential calculus. To this end, we assume smooth extensions of a scalar function $p:\, \Gamma \to \mathbb{R}$ or a vector function $\bu:\, \Gamma \to \mathbb{R}^3$ to a neighborhood $\mathcal{O}(\Gamma)$  of $\Gamma$. For example, one can consider extension along the normal directions using the closest point mapping $\bp:\,\mathcal{O}(\Gamma)\to \Gamma$. The surface gradient and covariant derivatives on $\Gamma$ are then defined as $\nablaG p=\bP\nabla p^e$ and  $\nabla_\Gamma \bu\coloneqq \bP \nabla \bu^e \bP$. The definitions  of surface gradient and covariant derivatives are  independent of a particular smooth extension of $p$ and $\bu$ off $\Gamma$.
The surface rate-of-strain tensor \cite{GurtinMurdoch75} is given by
 $E_s(\bu)\coloneqq\frac12(\nabla_\Gamma \bu + \nabla_\Gamma \bu^T)$,
and  the surface divergence operators for a vector $\bu: \Gamma \to \R^3$ and
a tensor $\bA: \Gamma \to \mathbb{R}^{3\times 3}$ are defined as
\[
 \divG \bu  \coloneqq \tr (\gradG \bu), \qquad
 \divG \bA  \coloneqq \left( \divG (\be_1^T \bA),\,
               \divG (\be_2^T \bA),\,
               \divG (\be_3^T \bA)\right)^T,
               \]
with $\be_i$ the $i$th basis vector in $\R^3$.

The conservation of momentum for a thin material layer together with the inextensibility condition and assumption of vanishing normal motions (geometric equilibrium)  leads to  the following surface Navier--Stokes problem:
Given area forces $\mathbf{f} \in L^2(\Gamma)^3$, with $\mathbf{f}\cdot\bn=0$, find a vector field $\bu:\, \lbrack0, T\rbrack\times\Gamma \to \R^3$, with $\bu(0, \cdot) = \bu_0$, $\bu\cdot\bn =0$, and $p:\, \lbrack0, T\rbrack\times\Gamma \to \R$ such that
\begin{equation}  \label{NSE}
\left\{
\begin{split}
  \rho\left(\frac{\partial\bu}{\partial t}+(\nabla_\Gamma \bu)\bu\right) - 2\nu\bP \divG (E_s(\bu))+\nabla_\Gamma p &=  \blf \quad \text{on}~\Gamma,   \\
  \divG \bu & =0 \quad \text{on}~\Gamma.
  \end{split}\right.
\end{equation}
Here $\bu$ is the tangential fluid velocity, $p$ is the surface fluid pressure, $\rho$ and $\nu$ are density and viscosity coefficients. We further assume $\nu$, $p$ and $\blf$ are re-scaled so that $\rho=1$.

\rev{
For time discretization, we assume a constant time step $\Delta t={T}/{N}$ and adopt the notation $\bu^k(\vect x)$ for velocity solution at time $t^k=k\Delta t$, and similar for $p^k$. A semi-implicit time-stepping scheme for~\eqref{NSE} reads: Given $\bu^{k-1}$, $\bu^{k-2}$, find $\bu^k$ s.t. $\bu^k\cdot\bn =0$ and $p^k$ solving
\begin{equation} \label{semi-discrete}
\left\{
\begin{split}
  \left[\bu\right]_t^{k}+(\nabla_\Gamma \bu^k)\widetilde\bu^{k} - 2\nu\bP \divG (E_s(\bu^k))+\nabla_\Gamma p^k &=  \blf^k \quad \text{on}~\Gamma,   \\
  \divG \bu^k & =0 \quad\:\: \text{on}~\Gamma
  \end{split}\right.
\end{equation}
for $k=2,3,\dots,N$.
In numerical experiments we employ the  second order method with
\begin{equation}\label{BDF2}
\left[\bu\right]_t^{k} =\frac{3\bu^{k}-4\bu^{k-1}+\bu^{k-2}}{2\Delta t},\quad \widetilde\bu^{k}=2\bu^{k-1}- \bu^{k-2},
\end{equation}
but this particular choice has little effect on the properties of the resulting linear systems.}

We see that on each time step the linearized problem in \eqref{semi-discrete} is the Oseen-type system
\begin{align}\begin{split}\label{split:tang:oseen}
	\alpha\,\vect u + (\nabla_\Gamma\vect u)\bw - 2\nu\,\vect P\divG( E_s(\vect u)) + \nabla_\Gamma p &= \vect f,\\
	\divG \vect u &= 0
\end{split}\end{align}
with $\alpha \simeq \Delta t^{-1}$,  the wind field~$\vect w:=\widetilde\bu^{k}$, and  $\vect f$ collects contributions
from the area forces and from the previous step velocities of the
discretizated time derivative. Hence the resulting system of linear algebraic equations is non-symmetric and of saddle-point type with properties resembling those of the planar Oseen system, well-studied problem; see, e.g. \cite{benzi2005numerical,ElmanBook}. In particular, the problem is increasingly hard to solve when $\nu$ goes to zero.
\rev{One way to avoid this increasing complexity of the linear algebra system is to lag the entire inertia term in time, e.g. to replace the second term in \eqref{semi-discrete} with $(\nabla_\Gamma \widetilde\bu^k)\widetilde\bu^{k}$, ending up with a symmetric Stokes-type problem, same on each time step. However, numerical stability of such implicit-explicit scheme is known~\cite{Temam84} to impose a time step restriction of the form~ $\Delta t\le c(\nu) h^d$, where $d=2$ for two-dimension flows and $c(\nu)$ decreasing for  $\nu\to0$. This leads to a serious growth of computational costs  for small $h$ and $\nu$ despite the ease of linear algebra. In contrast, \eqref{semi-discrete} is unconditionally stable~\cite{GR} (approximation analysis suggest $\Delta t\approx h$),  and our strategy here is to alleviate computational  costs by employing a more sophisticated linear algebra solve re-enforced by the recycling algorithm.}

A weak formulation of \eqref{split:tang:oseen} requires the closed subspace of $H^1(\Gamma)^3$  consisting of \emph{tangential} vector fields,
$ 
\bV_T \coloneqq \left\lbrace \bv \in H^1(\Gamma)^3 \mid \bv \cdot \bn =0 \right\rbrace.
$ 
$\bV_T$ serves as the velocity trial and  test space in a  weak formulation.
\rev{If desired, the tangentiality constraint can be relaxed in a penalty weak formulation~\cite{Jankuhn1}} that allows for the larger velocity space: $\{\bv\in L^2(\Gamma)^3\,:\, \bP\bv\in\bV_T\}$.
The well-posedness of both formulations relies on the surface Korn inequality (see (4.8) in \cite{Jankuhn1}):
 There exists a constant $c_K > 0$ such that
\begin{equation} \label{Korn}
\Vert \bu \Vert_{L^2(\Gamma)} +\Vert E_s(\bu) \Vert_{L^2(\Gamma)} \geq c_K \Vert \bu \Vert_{H^1(\Gamma)} \qquad \text{for all } \bu \in \bV_T,
\end{equation}
and the  inf-sup condition (Lemma 4.2 in \cite{Jankuhn1}):
There exists a constant $c_0>0$ such that the following holds:
\begin{equation}\label{LBB}
\inf_{p \in L^2_0(\Gamma)} \sup_{\bv \in \bV_T} \frac{\int_\Gamma q \divG \bu \, ds}{\Vert \bv\Vert_{H^1(\Gamma)} \Vert p \Vert_{L^2(\Gamma)}} \geq c_0.
\end{equation}

For the discretization of \rev{\eqref{split:tang:oseen}} we apply the  trace $\vect P_2$--$P_1$ FEM~\cite{olshanskiiinf}.
To apply the method, assume $\Gamma$ is strictly contained in a polygonal domain $\Omega \subset \R^3$ and consider a family $\{\T_h\}_{h >0}$ of shape regular tetrahedral tessellations  of $\Omega$.
Tetrahedra that have a nonzero intersection with $\Gamma$ are collected in the set
denoted by $\T_h^\Gamma$ with the characteristic mesh size  $h$, and  $\overline{\OGamma} = \bigcup\{\overline{T}\,:\,T\in \T_h^\Gamma \}$.
On $\overline{\OGamma}$ we consider the standard $H^1$-conforming finite element spaces of degree $k$,
\[
V_h^k=\{v\in C(\OGamma)\,:\, v\in P_k(T)~\text{for any}~T\in\T_h^{\Gamma}\},\quad \text{with}~k=1,2.
\]
The velocity and pressure bulk finite element spaces are then defined to be
\[
\bV_h \coloneqq (V_h^2)^3, \quad Q_h \coloneqq V_h^1\cap L^0_2(\Gamma).
\]
The finite element formulation uses the restrictions (traces) of these spaces on $\Gamma$. Note that  traces of vector functions from $\bV_h$ does not necessarily satisfy the tangentiality  $\bu\cdot\bn=0$ condition.
It is not straightforward, if possible at all, to build an $H^1$-conforming finite element method which is also conformal with respect to the tangentiality condition. Therefore,  the tangentiality condition will be enforced weakly by the penalty method.
We shall also need an extension  of the normal vector from $\Gamma$ to $\OGamma$, which we define as $\bn= \nabla d$, where $d$ is the signed distance function to $\Gamma$. In practice, $d$ is often not available and an approximation is used.
We introduce the following finite element bilinear forms:
\begin{align*}
 a(\bu,\bv) &\coloneqq
 	\rev{\int_{\Gamma}\left(\alpha \bu\cdot\bv+\bv^T(\nabla_\Gamma\vect u)\bw\right)\,ds} + 2\nu\int_\Gamma  E_s( \bu):  E_s( \bv)\, ds\\
 	&+ \tau \int_{\Gamma}(\bn\cdot\bu) (\bn\cdot\bv) \,ds + \rho_u \int_{\OGamma} [(\bn\cdot\nabla) \bu] \cdot [(\bn\cdot\nabla) \bv]  \, dx,\\
 \rev{b(\bv,q)} &\coloneqq \rev{-\int_{\Gamma} q \divG \bv \, ds = \int_{\Gamma} \nablaG q \cdot \bv \, ds}, \\
 \gamma(\bu,\bv)&\coloneqq \widehat{\gamma}\,\int_\Gamma \divG\vect u \divG\vect v\,ds, \quad s(p,q) \coloneqq \rho_p  \int_{\OGamma} (\bn\cdot\nabla p)(\bn\cdot\nabla q) \, dx.
\end{align*}
The forms are well defined for $p,q \in H^1(\OGamma)$, $\bu,\bv \in H^1(\OGamma)^3$. \rev{The finite element formulation of \eqref{split:tang:oseen} then reads: Find $\{\bu_h, p_h\} \in \bV_h \times Q_h$ solving
\begin{equation} \label{discrete}
 \begin{aligned}
	a(\bu_h,\bv_h) + \gamma(\rev{\bu_h, \bv_h}) + b(\bv_h,p_h) & =(\blf,\bv_h),  \\
	b(\bu_h,q_h)- s(p_h,q_h) & = 0
 \end{aligned}
\end{equation}
for all  $\bv_h \in \bV_h$ and $q_h \in Q_h$.}
In the definition of the forms,  $\tau>0$ is a penalty parameter and $\rho_u\ge0$, $\rho_p\ge0$ are stabilization parameters, which we set according to~\cite{olshanskiiinf}:
\begin{equation} \label{param}
\tau=h^{-2},\quad \rho_p=h, \quad \rho_u=1.
\end{equation}
Note that both stabilization terms vanish if $\bu$ and $p$ is the surface solution extended along normal directions to a neighborhood of $\Gamma$. This makes the finite element formulation consistent. At the same time, both terms add to the finite element formulation  additional ``stiffness'' in the normal direction. This
allows to eliminate the dependence of the resulting algebraic systems condition number on the position of $\Gamma$ in the background mesh, the idea first suggested for trace FEM in \cite{burman2015stabilized} and later explored in many publications on unfitted FEM for surface PDEs, e.g.~\cite{burman2017cut,burman2018cut,grande2017higher} (the particular choice of stabilization terms varies in the literature). We refer to  \cite{olshanskiiinf} for the further discussion of the role of these terms in the context of unfitted $\vect P_2$--$P_1$ elements, \rev{the proof of the well-posedness of the finite element formulation, error analysis,}  and the proof of $\Gamma$-independent estimates on the condition numbers of the velocity and pressure matricies.

The fourth term in \eqref{discrete} is the surface analogue of the grad-div stabilization~\cite{olshanskii2002low,olshanskii2004grad}. We further set 
$ 
\widehat{\gamma}=1
$ 
 unless it is stated otherwise and write $\gamma = \widehat{\gamma}$ to simplify the notation. \rev{We do not study the dependence of optimal $\gamma$ on other parameters of the finite element formulation. It is known~\cite{olshanskii2009grad} that there is a wide range of quasi-optimal $\gamma$-s, where the solution quality is almost insensitive to the variation of the parameter. Hence $\gamma$ can be taken smaller or larger depending on other considerations. For simplicity we adopt $\gamma=1$ for the full AL approach and mesh-dependent $\gamma$ for the modified AL approach; see the next section.}
It is interesting that no other stabilization was found to be necessary for computations with high Reynolds numbers.
A possible explanation is that tangential flows do not produce boundary layers (on a closed surface) and, in addition, the grad-div term by itself is known to dissipate excessive energy  in under-resolved simulation~\cite{olshanskii2009grad,john2010numerical}.

We conclude this section noting that the  implementation requires the integration of polynomial functions over $\Gamma$. In practice this is avoided by approximating $\Gamma$ by some $\Gamma_h$ which admits exact quadrature rules. The quantification of the introduced geometric inconsistency in the case of the Stokes problem and $\vect P_k$--$P_{k-1}$, $k=2,3,\dots$, trace elements is given in \cite{jankuhn2020error}.

\section{System of linear algebraic equations and preconditioning}\label{s:matr}
We now turn to the matrix form of the discretized surface Oseen system and define the velocity, pressure stabilization and divergence constraint matrices:
\[
\bA_{i,j}= a(\psi_j,\psi_i)+ \gamma(\rev{\psi_j,\psi_i}),\qquad
\bC_{k,m}=s(\xi_k,\xi_m),\qquad \bB_{k,i}= b(\psi_i,\xi_k),
\]
where $\{\psi_i\}$ and $\{\xi_k\}$ are the velocity and pressure nodal basis  functions spanning $\bV_h$ and $Q_h$, respectively.
After arranging velocity degrees of freedom first and pressure degrees of freedom next, we arrive at the system with the $2\times2$-block matrix:
\begin{equation}\label{LAS}
\underbrace{
\left[\begin{matrix}
\bA & \bB^T \\
\bB & -\bC
\end{matrix}\right]}_{\mathbf{ \mathcal{A}}}
\left[
\begin{matrix}
\mathrm{u}\\
\mathrm{p}
\end{matrix}\right]= \left[ \begin{matrix}
\mathrm{f}\\
\mathrm{0}
\end{matrix}\right].
\end{equation}
An important matrix related to the above system is the pressure  Schur complement
$ \bS=\bB \bA^{-1} \bB^T + \bC$,
which results after elimination of the velocity unknowns  from the system. A preconditioner for $\bS$ is a necessary ingredient for most iterative solvers that exploit the block structure of $\mathcal{A}$. Following a common practice~\cite{ElmanBook} we consider the block-triangle right preconditioner for $\mathcal{A}$:
\begin{equation}\label{split:tang:prec}
	\mathcal P \coloneqq \begin{pmatrix}
		\widehat{\vect A} & \vect B^T \\
		& \widehat{\vect S}\phantom{^T}
	\end{pmatrix},
\end{equation}
where $\widehat{\vect A}$ and $\widehat{\vect S}$ are preconditioners for $\vect A$ and ${\vect S}$, respectively.

For the next step, we define the surface pressure mass matrix $\bM_p$ and the pressure  Laplace--Beltrami matrix~$\vect L_p$:
\[
(\bM_p)_{k,m}=\int_\Gamma\xi_k\xi_m\,ds, \qquad (\bL_p)_{k,m}=\int_\Gamma\nablaG \xi_k\cdot \nablaG \xi_m\,ds.
\]
For the surface Stokes problem ($\alpha=0$, $\bw=0$, $\gamma=0$, $\nu=1$) matrix $\bS$ is spectrally equivalent to the stabilized pressure mass matrix $\bM_p+\bC$; see~\cite{olshanskiiinf}. However, for $\bw\neq0$, $\gamma=0$, and $\nu\to0$ the problem of building a suitable preconditioner for $\bS$ is known to be particular difficult. To circumvent it, the authors of \cite{benzi2006augmented} introduced an augmentation to the $\bA$ block of the system \rev{replacing $\bA$ with $\bA+\gamma\bB^T \vect M_p^{-1}\bB$. Such augmentation is not algebraically consistent in our case, since $\bC\neq0$ and so $\bB\,\mathrm{u}\neq0$. We note that $\bC\neq0$ is a typical situation for many unfitted \textit{inf-sup stable} FEM discretizations of the (Navier--)Stokes equations (both in volumes and on surfaces) as well as for stabilized elements~\cite{benzi2011analysis}. Hence, we suggest to introduce the augmentation on the finite element level, i.e. to add the grad-div term.}

For the planar Oseen problem discretized with standard $\vect P_2$--$P_1$ elements one can show that the   Schur complement of the algebraically augmented matrix is spectrally equivalent to the pressure mass matrix scaled by  $(\gamma+\nu)^{-1}$ for sufficiently large $\gamma$~\cite{benzi2006augmented}.
\rev{
We show here that similar property holds for the trace FEM and when the algebraic augmentation is replaced by the grad-div stabilization
and so the augmentation term is not of the $\gamma\bB^T\vect M_p^{-1}\bB$ form. More precisely, assume $\alpha=0$, $\nu\le1$, $\|\bw\|_{L^\infty(\Gamma)}=1$ and the skew-symmetric discretization of the advection term \eqref{aux691} (these assumptions are standard for the analysis, but can be relaxed for the expense of extra technical details), then the eigenvalues of
\begin{equation}\label{eig_pr}
 \vect S\,\mathrm p =\lambda \big((\nu + \gamma)^{-1}\,\vect M_p + \vect C\big)\,\mathrm p
\end{equation}
satisfy the following bound
\begin{equation}\label{bound}
  c\frac{\gamma+\nu}{\gamma+\nu^{-1}}\le\Re(\lambda),\quad |\lambda| \le \frac{\nu+\gamma}{\gamma},
\end{equation}
with some positive $c$ independent of problem parameters and the position of $\Gamma$ in the background mesh. $\Re(\lambda)$ above denotes the real part of $\lambda$. We include the proof in the Appendix. We see that for large enough $\gamma$ all eigenvalues are contained in a box in the right half-plane with the bounds independent of $\nu$.
Motivated by \eqref{bound} we} define  the Schur complement preconditioner through its inverse as follows:
	\begin{equation}\label{split:tang:prec:schur}
	\widehat{\vect S}^{-1} \coloneqq \big((\nu + \gamma)^{-1}\,\vect M_p + \vect C\big)^{-1} + \big(\alpha^{-1}\,\vect L_p + \vect C\big)^{-1}
	\end{equation}
	 The second term is included to deal with the dominating reaction term in the Oseen problem~\eqref{split:tang:oseen} if $\alpha\gg1$. This part of~$\widehat{\vect S}^{-1}$ resembles the Cahouet--Chabard preconditioner~\cite{cahouet1988some}. We apply several CG iterations to compute the action of $\big(\alpha^{-1}\,\vect L_p + \vect C\big)^{-1}$ and $\big((\nu + \gamma)^{-1}\,\vect M_p + \vect C\big)^{-1}$ on a vector. \rev{Alternatively, these matrices can be also one time factorized. Since the number of pressure degrees of freedom is much smaller than velocity ones, either choice marginally affects the total timings. Note also that $\alpha^{-1}\,\vect L_p + \vect C$ has a one-dimensional kernel, i.e. the subspace of
	 	constant pressures, which we easily handle by iterating in a proper subspace. Strictly speaking $\big(\alpha^{-1}\,\vect L_p + \vect C\big)^{-1}$ is the pseudo-inverse in our case. }

The augmentation has the downside of adding to the $(1, 1)$-block the term with a large nullspace. For larger $\gamma$ this makes the matrix $\bA$ poor conditioned and hinders the efficiency of standard iterative methods to evaluate ${\vect A}^{-1}$.
As a more flexible alternative we explore here direct LU factorization of $\bA$ (or its sub-blocks) and the reuse of the factors for several time steps. In pursuing this line, we consider two strategies of building $\widehat{\vect A}$:
\begin{enumerate}
  \item LU factorization of the full velocity block $\bA$ (full AL approach);
  \item Velocity unknowns are enumerated componentwise so that $\vect A$ attains the $3\times3$-block form.   $\widehat{\vect A}$ is obtained from ${\vect A}_{tr}$, the block upper-triangle part of $\vect A$, by applying LU factorization to each individual diagonal block of ${\vect A}_{tr}$. This corresponds to modified AL approach from~\cite{benzi2011modified}.
\end{enumerate}
The modified AL approach allows to factorize smaller and better structured matrices that have the structure of a stiffness matrix of a conforming FEM applied to an elliptic scalar PDE.
This enhanced efficiency comes with a price of slight $\nu$ and $h$ dependence of the preconditioner performance~\cite{benzi2011modified}. We shall see below that in the case of time-dependent 2D flow the price is very tolerable.

Same approaches, of course, apply to reusing ILU factorizations, but we fix our idea and consider below exact LU.  The surface fluid problem is essentially 2D and the number of velocity unknowns allows applying LU factorization. Furthermore, in a curvilinear metric the viscosity term does not simplify to Laplace operators for each velocity component, \rev{i.e. for tangential divergence free $\bu$ we note that in general $\bP \divG (E_s(\bu))\neq \Delta_\Gamma \bu$ with a componentwise Laplace--Beltrami operator $\Delta_\Gamma$.  Therefore, $\vect A$ does not have a block-diagonal structure for $\gamma=0$ and so} adding the $\gamma$-term does not \rev{change} the sparsity pattern of the matrix (in contrast to the augmentation in the planar case).

To make the algorithm precise, denote by $\bL(k)$ and $\bU(k)$ the LU factors of $\bA(k)$ at step $k$ of~\eqref{semi-discrete}. We let $\widehat{\vect A}=\bL(k)\bU(k)$ to be the preconditioner for all $\bA(k+m)$, $m=0,\dots,M$, where $M\ge0$ is the largest index such that
\begin{equation}\label{Recycle}
  \frac{\#\mbox{Iter}_{\rm FGMRES}(k+m)}{\#\mbox{Iter}_{\rm FGMRES}(k)}\le \kappa,\quad \text{for}~m=0,\dots,M,
\end{equation}
where $\kappa\ge1$ is a maximum allowed increase  of the preconditioned  FGMRES iterations without updating the preconditioner.

\section{Numerical simulations}\label{s_num}
We apply the trace FEM as described in section~\ref{sectioncont} to simulate  the mixing layer of  isothermal incompressible viscous surface fluid at several Reynolds numbers.
The setup resembles the  classical problem of the Kelvin--Helmholtz instability: For a detailed discussion of the  planar analogue we refer to~\cite{schroeder2019reference} and references therein. At higher Reynolds numbers the flow exhibits sharp internal layers and intensive vortical dynamics, offering a good test problem for both discretizations and flow solvers.

For discretization,  an initial triangulation $\T_{h_0}$  was build by  dividing $\Omega=(-\frac53,\frac53)^3$ into $2^3$ cubes and further splitting each cube into 6 tetrahedra with  $h_0 = \frac53$. Further, the mesh is  refined only close to the surface, and $\ell \in \mathbb{N}$ denotes the level of refinement so that $h_\ell = \frac53\,2^{-\ell}$. The  trace $\vect P_2$--$P_1$ Taylor--Hood finite element method with BDF2 time stepping is applied as described in \eqref{semi-discrete}--\eqref{BDF2}, with the choice of parameters in \eqref{param}. \rev{The choice of $\gamma$ is $\gamma=1$ in the full AL and a mesh-dependent $\gamma$ in the modified AL.}
No further stabilizing terms, e.g., of streamline diffusion type, were included in the method, since the computed solutions do not reveal any spurious modes. We would like the discretization error, which results from the approximation of $\Gamma$, to be consistent with the higher order interpolation properties of $\vect P_2$ elements~\cite{jankuhn2020error}.  To address this, we apply additional refinement to define a piecewise-linear surface $\Gamma_h$ used for the purpose of numerical integration as described in~\cite[section~6.3]{olshanskiiinf}.
 Software packages DROPS~\cite{DROPS} and Belos, Amesos from Trilinos~\cite{trilinos-website} were used for matrices assembling and algebraic solver execution, respectively. Because of the additional refinement used to define numerical quadratures, the matrix assembling time grows superlinear in our examples. The optimal complexity here can be obtained by using isoparametric higher order trace elements~\cite{lehrenfeld2016high}, not however implemented in the software we use.

\subsection{The Kelvin--Helmholtz instability problem setup}
There are very limited numerical studies of the Kelvin--Helmholtz (KH) problem for surface fluids. Examples of an isothermal  KH flows on  cylinder and on the unit sphere are given in~\cite{lederer2020divergence,jankuhn2020error}. Here we use the sphere example and for the first time simulate  the KH flow on a torus of revolution.

The design of numerical experiment for the sphere follows~\cite{lederer2020divergence,jankuhn2020error}.  For~$\Gamma=S^2$, let $\xi$ and $\zeta$ to be {renormalized} azimuthal and polar coordinates, respectively: $-1/2 \le \xi, \zeta < 1/2$. The tangent basis direction are~$\vect e_\xi \coloneqq \nabla_{\Gamma}\xi/\|\nabla_{\Gamma}\xi\|$ and $\vect e_\zeta \coloneqq \nabla_{\Gamma}\zeta/\|\nabla_{\Gamma}\zeta\|$. The initial velocity field is given by the counter-rotating upper and lower hemispheres with velocity speed approximately equal 1 closer to equator and vanishing at poles. The velocity field has a sharp transition  layer along equator, where we add perturbation to trigger the development of the vortical strip:
\begin{align}\begin{split} \label{kh:ini}
	{\vect u}_0(\xi, \zeta) &\coloneqq d(\zeta)(\tanh(2\,\zeta/\delta_0)\,\vect e_\xi + c_n\vCurl_\Gamma\psi), \\
	\psi(\xi, \zeta) &\coloneqq e^{-(\zeta/\delta_0)^2}\,\big(a_a\cos(m_a\,\pi\,\xi) + a_b\cos(m_b\,\pi\,\zeta)\big),
\end{split}\end{align}
where $d$ is the distance from $\Gamma$ to the $z$-axis. We take $\delta_0 \coloneqq 0.05$ (for $|z| \gtrsim \delta_0$ the velocity field is close to a rigid body rotation around the $z$-axis), $c_n \coloneqq 10^{-2}$ (perturbation parameter), and $a_a = 1$, $m_a = 16$, $a_b = 0.1$, $m_b = 20$ (perturbation magnitudes and frequencies). Note that ${\vect u}_0$ is tangential by construction. 
The Reynolds number $\text{Re} \simeq \nu^{-1}\delta_0$ is based on $\|\vect u\|_{\LInfSpace} \simeq 1$ and the initial layer width. We ran numerical simulations with $\nu = \frac12 10^{-k}$, for $k=3,4,5$, which corresponds to   $\text{Re}=10^{k}$, $k=2,3,4$.

\begin{figure}[ht!]
	\centering
		\includegraphicsw[.2]{{kh_t=0.cropped}.png}~
		\includegraphicsw[.2]{{kh_t=2.5.cropped}.png}~
		\includegraphicsw[.2]{{kh_t=5.cropped}.png}~
		\includegraphicsw[.2]{{kh_t=6.25.cropped}.png}
\vskip1ex
		\includegraphicsw[.2]{{kh_t=10.cropped}.png}~
		\includegraphicsw[.2]{{kh_t=12.5.cropped}.png}~
		\includegraphicsw[.2]{{kh_t=15.cropped}.png}~
		\includegraphicsw[.2]{{kh_t=20.cropped}.png}
\vskip1ex
\includegraphicsw[.4]{{sphere_kh_legend.cropped}.png}
	\caption{KH flow at $\text{Re}=10^{4}$. Snapshots of surface vorticity~${w_h = \Curl_{\Gamma_h}\vect u_h}$ for~${t \in \{ 0, 2.5, 5, 6.25, 10, 12.5, 15, 20 \}}$.  Click \href{https://youtu.be/EwF3vCgFhuI}{\underline{here}} for the full animation.}
	\label{fig:kh:curl}		
\end{figure}

In Figure~\ref{fig:kh:curl} we show several snapshots of the surface vorticity distributions starting from the initial condition. The solution reproduces the well known flow pattern of the planar KH instability development, which includes  the initial vortices formation in the layer followed by pairing and self-organization into larger vortices.  At the final simulation time we see two large  counter-rotating vortices. The dissipation rate for this solution was studied in~\cite{jankuhn2020error} confirming the qualitatively correct behaviour.

The second example we consider is KH instability on two-dimensional torus $\Gamma=\mathbb{T}^2 = \{ \bx\in\Omega \mid r^2 = x_3^2 + (\sqrt{x_1^2 + x_2^2} - R)^2\}$, with  $R= 1$ and $r= 0.5$. The coordinate system is  $(\rho, \phi,%
\theta)$, with
\[
    \bx = R\begin{pmatrix}\cos\phi \\ \sin\phi \\ 0 \end{pmatrix}
    + \rho\begin{pmatrix}\cos\phi\cos\theta \\ \sin\phi\cos\theta \\
    \sin\theta \end{pmatrix},
\]
where $\rho$-direction is normal to $\Gamma$, $\frac{\partial \bx}{\partial\rho}\perp\Gamma$ for $\bx\in\Gamma$.
In the torus coordinates, the initial velocity field is given by the same formula \eqref{kh:ini} with $\xi=\psi/(2\pi)$,
$\zeta=\theta/(2\pi)$, and $d(\zeta) \coloneqq d(\bx(\xi,\zeta))=\sqrt{x^2 + y^2} - 0.5$, $\bx=(x,y,z)$, so that $d(\zeta)$ vanishes on  the inner ring of the torus.

\begin{figure}[ht!]
	\centering
		\includegraphicsw[.29]{{img/tor_kh_t=0.cropped}.png}~
		\includegraphicsw[.29]{{img/tor_kh_t=2.3.cropped}.png}~
		\includegraphicsw[.29]{{img/tor_kh_t=8.4.cropped}.png}
\vskip1ex
		\includegraphicsw[.29]{{img/tor_kh_t=9.3.cropped}.png}~
		\includegraphicsw[.29]{{img/tor_kh_t=14.cropped}.png}~
		\includegraphicsw[.29]{{img/tor_kh_t=18.7.cropped}.png}
\vskip1ex
		\includegraphicsw[.29]{{img/tor_kh_t=25.7.cropped}.png}~
		\includegraphicsw[.29]{{img/tor_kh_t=30.4.cropped}.png}~
		\includegraphicsw[.29]{{img/tor_kh_t=35.cropped}.png}
\vskip1ex
\includegraphicsw[.4]{{tor_kh_legend.cropped}.png}
	\caption{Snapshots of surface vorticity~${w_h = \Curl_{\Gamma_h}\vect u_h}$ for KH on the $R=1,\, r=\frac12$ torus  for~${t \in \{ 0, 2.3, 8.4, 9.3, 14, 18.7, 25.7, 30.4, 35 \}}$.  Click \href{https://youtu.be/v0bmM-NqRqo}{\underline{here}} for the full animation.}
	\label{fig:kh_t:curl}		
\end{figure}

Figure~\ref{fig:kh_t:curl} visualizes the vorticity field of the KH flow on the torus for $Re=10^4$ ($\nu=\frac12\, 10^{-5}$). The initial stage of the vortical layer formation and small vortices pairing is similar to the case of the sphere and the plane. The different geometry (and topology) of the torus apparently affects the interaction of larger vortices. From the time of about $20$ units there are 4 large vortices formed, which further travel in the both toroidal and poloidal directions  without  pairing up to time $t=45$, \rev{after $t=45$ the motions loses any apparent axial symmetry and becomes rather complex (see the animation)}.

\subsection{Solver performance}

\begin{table}[ht]
	\centering
	\caption{Solver statistics for $\nu = \frac 12\times10^{-4}$, varying $h$ and $\alpha\rev{\in \{24,48,96\}}$. Full AL preconditioner; $\Gamma$ is the sphere. \rev{The total number of time steps for three mesh levels reported are 320, 640, and 1280, respectively.}}
	\label{tab:h}
	\begin{tabular}[1.8]{|c|c|c|c|c|c|c|c|c|}
		\hline
		\multirow{2}{*}{\# d.o.f.} & \multirow{2}{*}{$T_\text{asmbl}$} & \multirow{2}{*}{\% factor steps} & \multicolumn{3}{c|}{``fresh'' LU steps} & \multicolumn{3}{c|}{all steps} \\
		\cline{4-9}
		& & & $N_\text{iter}$ & $T_\text{factor}$ & $T_\text{linsol}$ & $N_\text{iter}$ & $T_\text{factor}$ & $T_\text{linsol}$ \\
		\hline
		$51526$	&	$3.55$	&	$3.44$	&	$9.00$	&	$1.72\times	10^1$	&	$0.53$	&	$33.53$	&	$0.59$	&	$2.15$	\\ \hline
$203998$	&	$18.6$	&	$1.88$	&	$8.33$	&	$1.76\times	10^2$	&	$2.55$	&	$32.77$	&	$3.30$	&	$10.3$	\\ \hline
$819862$	&	$180$	&	$0.86$	&	$7.55$	&	$1.73\times	10^3$	&	$12.1$	&	$29.86$	&	$14.9$	&	$49.9$	\\ \hline 
	\end{tabular}
\end{table}

\begin{figure}
	\centering
	\begin{subfigure}{.5\linewidth}
		\centering
		\includegraphicsw[0.9]{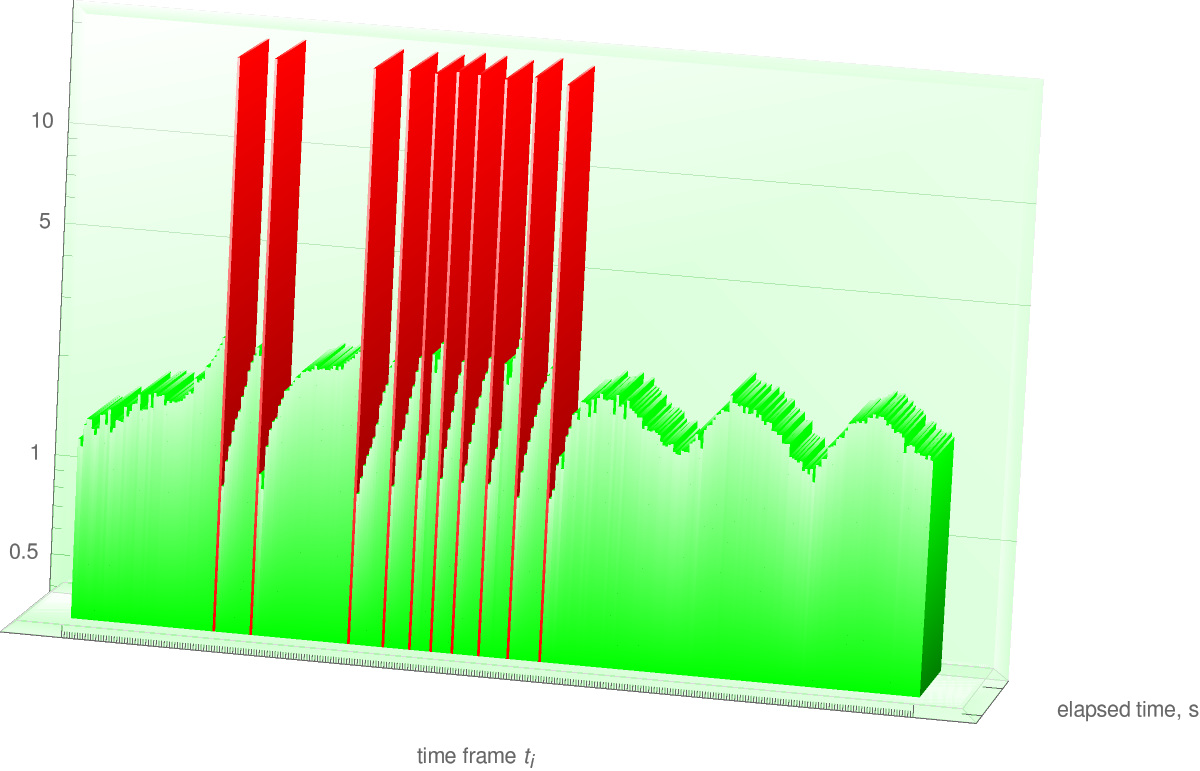}
		\caption{\# d.o.f. = 51526, $n_{LU} = 10$}
	\end{subfigure}%
	\begin{subfigure}{0.5\linewidth}
		\centering
		\includegraphicsw[0.9]{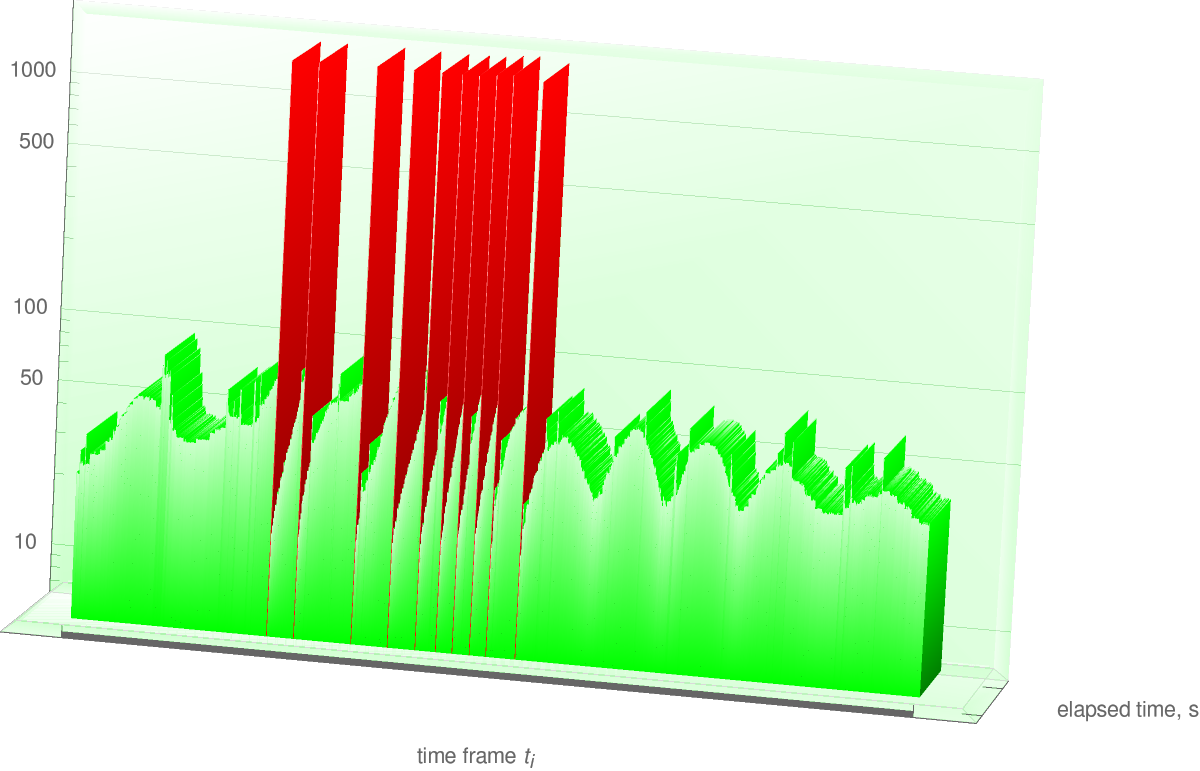}
		\caption{\# d.o.f. = 819862, $n_{LU} = 10$}
	\end{subfigure}
	\caption{Computation time seconds (factorization and linear solve) in log-scale vs. time $t_i$ for two different meshes. Red bars correspond to time steps for which new factors are computed; $n_{LU}$ is a number of  steps with factorizations. Full AL preconditioner; $\Gamma$ is the sphere.} \label{fig:h}
\end{figure}

In the first series of experiments we fixed the viscosity parameter equal $\frac1210^{-4}$ and vary the mesh size. We consider three levels of mesh refinement and the number of unknowns grows by  a factor of four from one level to the next one. Since $\alpha\simeq 1/\Delta t$ and $\Delta t\simeq h$, the parameter $\alpha$ in \eqref{split:tang:oseen} increases two times for each refinement level.  The FGMRES with full AL preconditioner was applied to solve the system of algebraic equations on each time step of \eqref{discrete}. We use zero vector as initial guess and the drop of residual by a factor of~$10^8$ as the stopping criterium.
Table~\ref{tab:h} summarizes the solver averaged statistics over the time of simulation~ $t\in[0,20]$.  We see that the percentage of re-initializations of the preconditioner (this is when we compute new LU factors) is small and decreases for finer mesh levels. The later can be due to the growth of $\alpha$ and because the the diffusion term plays more significant role for smaller $h$. The choice of $\kappa=5$ in \eqref{Recycle} keeps the average number of FGMRES iterations about 30 with very slight variation among refinement levels.  To compare, the number of FGMRES iterations with `fresh' LU factorization in $\widehat{\bA}$ is about 8. As expected, factorizing the matrix $\widehat{\bA}$ is by far the most computationally expensive procedure (cf. $T_\text{factor}$ in ``fresh LU steps'' table section). However, due to the heavy and efficient recycling, overall the expense of the factorization is minor compared to the  iterations cost (compare $T_\text{factor}$ and $T_\text{linsol}$  in ``all steps'' table section). This allows to keep the \textit{averaged} computation cost of the linear solver comparable and even less than the cost of matrix assembling.
This balance is further visualized in Figure~\ref{fig:h} for two mesh levels, where we see that the time steps with updated LU factors are more expensive but rare. It is interesting to note that most updates are needed for $t\in[4,10]$, when vortixes are paring. As we discussed above, the numerical integration that we use causes the assembling time to grow superlinear with respect to \#d.o.f.: This is the specific of the flow problem posed on a manifold and software we use for matrix assembly.

\begin{table}[ht]
	\centering
	\caption{Solver statistics for fixed $h$ (\# d.o.f. = 819862) and $\alpha\simeq100$, varying $\nu$.  Full AL preconditioner; $\Gamma$ is a sphere.}
	\label{tab:nu}
	\begin{tabular}[1.8]{|c|c|c|c|c|c|c|c|c|}
		\hline
		\multirow{2}{*}{$\nu$} & \multirow{2}{*}{$T_\text{asmbl}$} & \multirow{2}{*}{\% factor steps} & \multicolumn{3}{c|}{``fresh'' LU steps} & \multicolumn{3}{c|}{all steps} \\
		\cline{4-9}
		& & & $N_\text{iter}$ & $T_\text{factor}$ & $T_\text{linsol}$ & $N_\text{iter}$ & $T_\text{factor}$ & $T_\text{linsol}$ \\[0.4ex]
		\hline
		$\frac12\times	10^{-3}$	&	177	&	0.391	&	$7.00$	&	$1.59\times	10^3$	&	$9.87$	&	$32.87$	&	$6.23$	&	$50.7$	\\[0.5ex] \hline
$\frac12\times	10^{-4}$	&	180	&	0.859	&	$7.55$	&	$1.73\times	10^3$	&	$12.1$	&	$29.86$	&	$14.9$	&	$49.9$	\\[0.5ex] \hline
$\frac12\times	10^{-5}$	&	198	&	0.938	&	$7.75$	&	$1.81\times	10^3$	&	$12.4$	&	$31.95$	&	$17$	&	$50.2$	\\ \hline 
	\end{tabular}
\end{table}

\begin{figure}[H]
	\centering
	\begin{subfigure}{.5\linewidth}
		\centering
		\includegraphicsw[0.9]{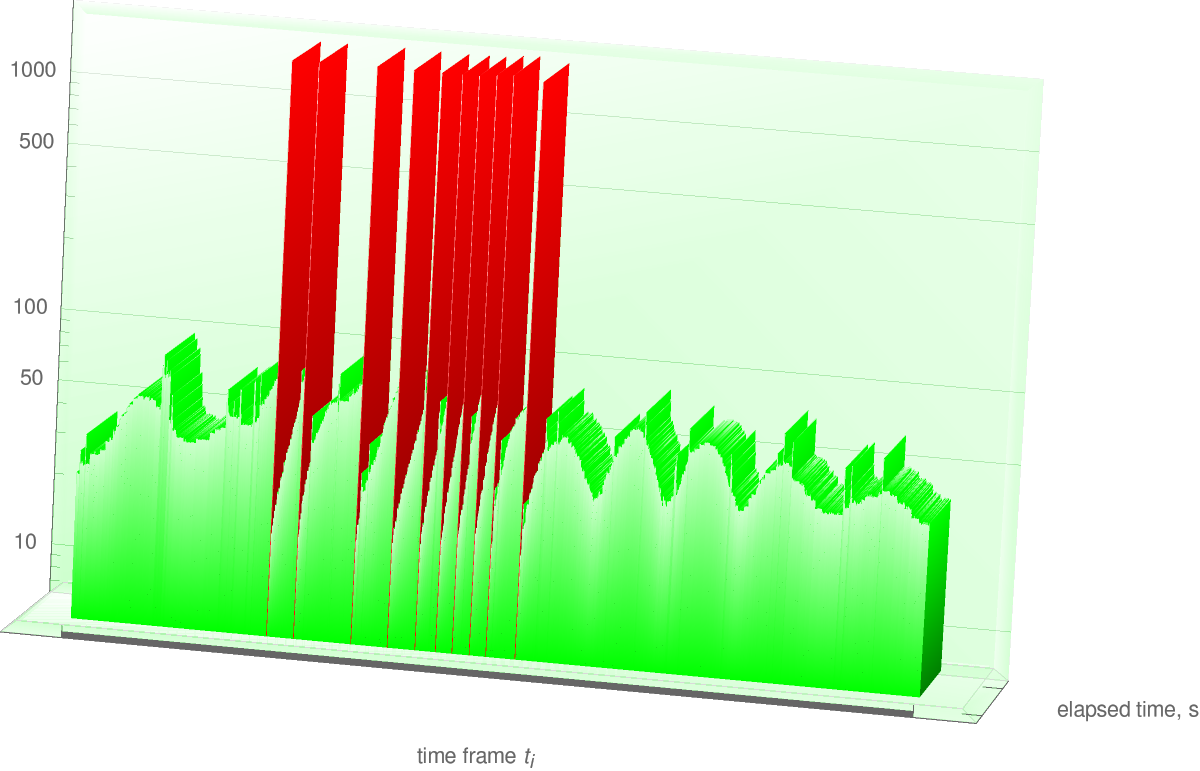}
		\caption{$\nu = \frac12\times10^{-4}$, $n_{LU} = 10$}
	\end{subfigure}%
	\begin{subfigure}{0.5\linewidth}
		\centering
		\includegraphicsw[0.9]{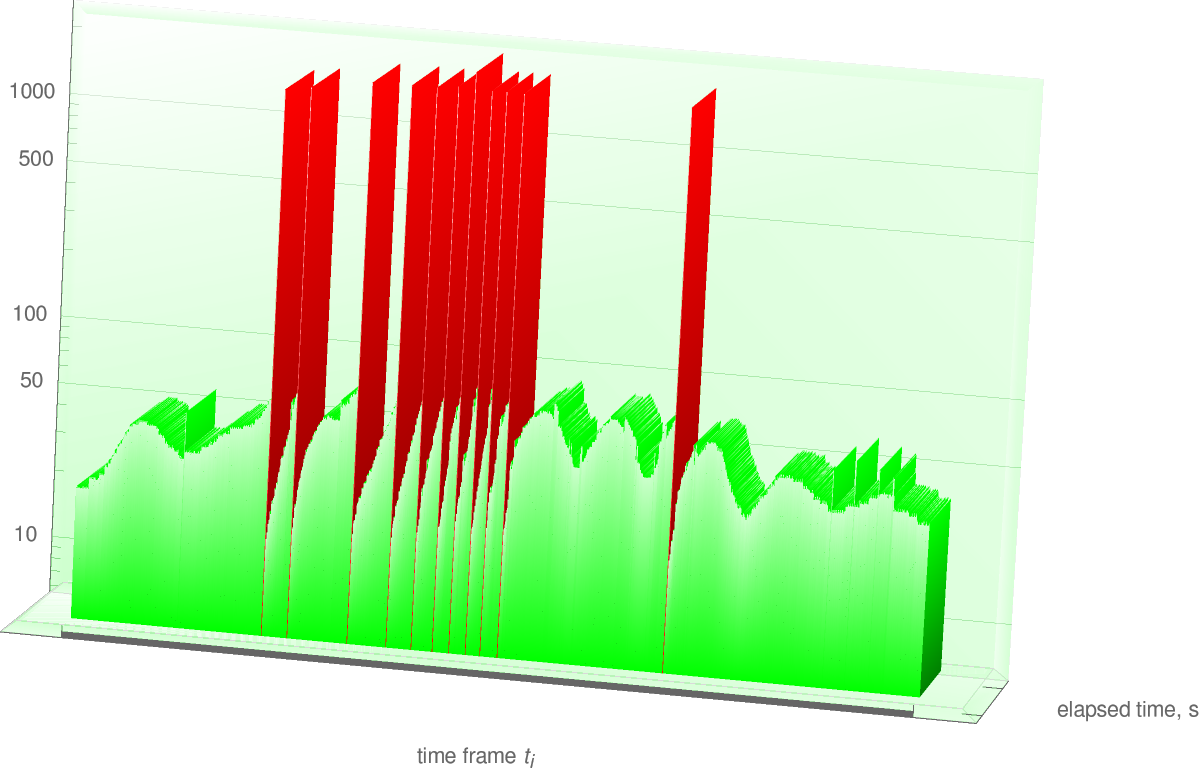}
		\caption{$\nu = \frac12\times10^{-5}$, $n_{LU} = 11$}
	\end{subfigure}
	\caption{Computation time seconds (matrix assembly and linear solve)  in log-scale  vs. time $t_i$ for two different viscosity parameter values. Red bars correspond to time steps for which factorization is performed; $n_{LU}$ is a number of factorization steps. Full AL preconditioner; $\Gamma$ is a sphere.
} \label{fig:nu}
\end{figure}

We repeat the simulation of the KH problem on the sphere but now for several values the viscosity parameter and the finest discretization level. All parameters of the algebraic solver are the same as above. The averaged statistics of the solver for this set of experiments are summarized in Table~\ref{tab:nu}. It appears that the solver is remarkably robust with respect to the viscosity parameter.  For higher Reynolds numbers we see only a slight increase of the percentage of time steps, where the preconditioner is updated by the new LU factors.
Figure~\ref{fig:nu} illustrates the balance between computationally expensive but rare steps with updated preconditioner and the majority of calculations with the recycled AL preconditioner.

\begin{table}[H]
	\centering\footnotesize
	\caption{ \label{tab:h:tor}
Solver statistics for fixed $\nu = \frac 12\times10^{-5}$, varying $h$ and $\alpha \rev{\in \{24,48,96\}}$.   Modified AL preconditioner; $\Gamma$ is a torus.
\rev{The total number of time steps for three mesh levels reported are 960, 1920, and 3840, respectively.} }
	\label{tabh:tor}
	\begin{tabular}[1.8]{|c|c|c|c|c|c|c|c|c|c|}
		\hline
		\multirow{2}{*}{\# d.o.f.} & \multirow{2}{*}{$\gamma$} & \multirow{2}{*}{$T_\text{asmbl}$} & \multirow{2}{*}{\% factor steps} & \multicolumn{3}{c|}{``fresh'' LU steps} & \multicolumn{3}{c|}{all steps} \\
		\cline{5-10}
		& & & & $N_\text{iter}$ & $T_\text{factor}$ & $T_\text{linsol}$ & $N_\text{iter}$ & $T_\text{factor}$ & $T_\text{linsol}$ \\
		\hline
		$78244$		&	$0.028$	&	$7.65$	&	$6.25\times	10^{-1}$	&	$37.33$	&	$1.63$	&	$2.94$	&	$68.31$	&	$1.02\times	10^{-2}$	&	$5.02$	\\ \hline
$315792$	&	$0.020$	&	$38.3$	&	$4.17\times	10^{-1}$	&	$42.75$	&	$18.1$	&	$16.3$	&	$75.32$	&	$7.55\times	10^{-2}$	&	$24.6$	\\ \hline
$1279180$	&	$0.014$	&	$324$	&	$3.65\times	10^{-1}$	&	$64.71$	&	$181$	&	$97.6$	&	$75.51$	&	$6.61\times	10^{-1}$	&	$112$	\\ \hline 
	\end{tabular}
\end{table}

\begin{figure}[H]
	\centering
	\begin{subfigure}{.5\linewidth}
		\centering
		\includegraphicsw[0.9]{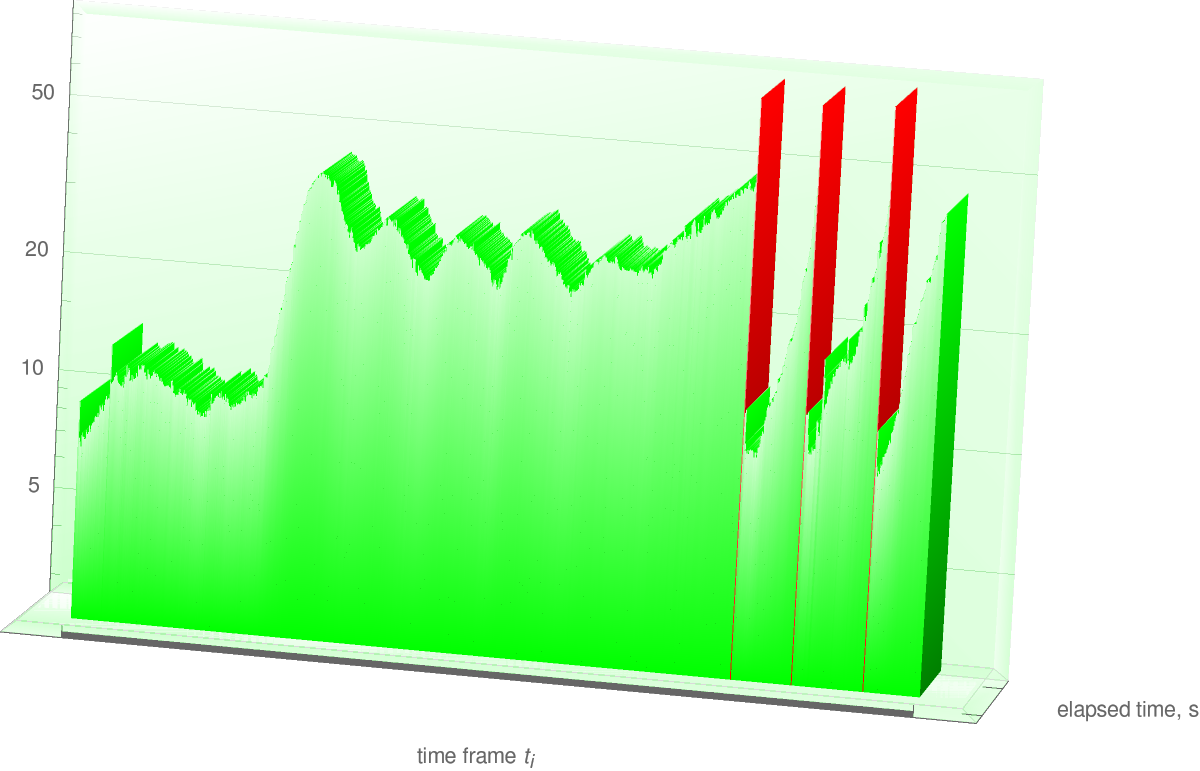}
		\caption{\# d.o.f. = 315792, $n_{LU} = 3$}
	\end{subfigure}%
	\begin{subfigure}{0.5\linewidth}
		\centering
		\includegraphicsw[0.9]{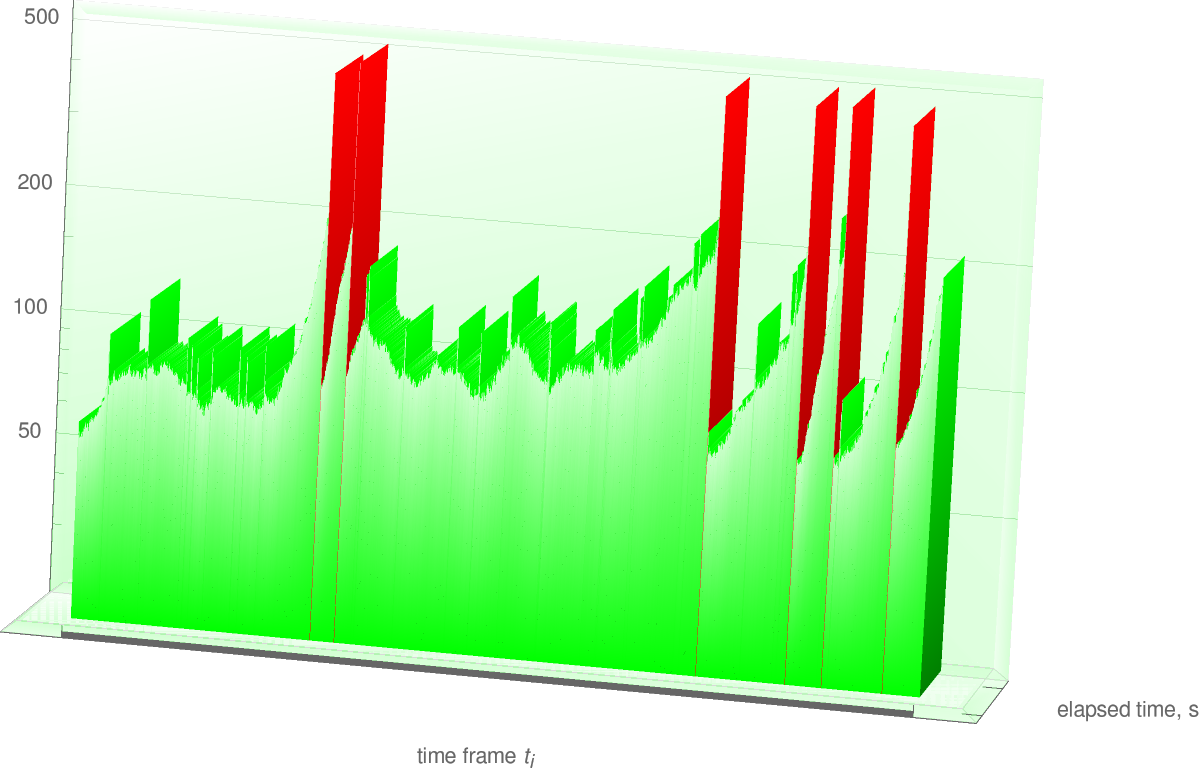}
		\caption{\# d.o.f. = 1279180, $n_{LU} = 6$}
	\end{subfigure}
	\caption{Computation time seconds (matrix assembly and linear solve) in log-scale vs. time $t_i$. Red bars correspond to time steps for which new factors are computed; $n_{LU}$ is a number of  steps with factorizations. Modified AL preconditioner; $\Gamma$ is a torus.} \label{fig:h:tor}
\end{figure}

We now consider the surface  KH flow on torus. For the given values of outer and inner radius the surface area of torus is approximately 1.57 times the surface of the unit sphere. This explains, why we get larger systems in terms of the number of degrees of freedom for the same levels of refinement in this example.  This makes the problem naturally suitable for testing the recycling strategy with modified AL preconditioner. In general, the  modified AL preconditioner is less robust with respect to $\nu$ and $h$, so its efficient use needs some tuning.
Following  recommendations in~\cite{benzi2011modified} we find optimal value for $\gamma$ on a coarse level and then apply $1/\sqrt{2}$-rule to scale it for finer mesh levels. This leads us to $\gamma_3 = 0.04$ for the third refinement level and $\gamma_i =2^{\frac{3-i}{2}}\gamma_3 $, $i=4,5,6$, for refinement levels from 4 to 6. These are the refinement levels we use to report the solver statistics in  Table~\ref{tab:h:tor}. In this experiment, we take the velocity field and pressure from the previous time step as the initial guess in FGMRES and relax the stopping criterium to the relative drop of residual by $10^6$, $\|\vect r_i\|_2 < 10^{-6} \|\vect b\|_2$. The number of iterations increased compared to the full AL preconditioner, since only the block upper-triangle part of the matrix~$\bA$ is used to define $\widehat{\bA}$. We also see a slight increase of the iteration number for $h$ getting smaller, which is also the observation in~\cite{benzi2011modified}. The overall computation time is dominated by the matrix assembly because of the non-optimal numerical quadrature, as discussed above. If we take the time of assembly off the table, then the recycling strategy turns out to be very effective also with the modified AL preconditioner. The average time of factorization per one solve is negligible and each factorization is more efficient in terms of time and memory requirements since it is done for each individual velocity block.
Figure~\ref{fig:h:tor} illustrates the balance between computationally expensive but rare steps with updated preconditioner and the majority of calculations with the recycled modified AL preconditioner.

\rev{It is out of scope for this paper to carry out a systematic comparison of the full and modified AL preconditioners. Results in this direction can be found in \cite{benzi2011modified}. For the factorize--recycle framework introduced here, our general recommendation is the following: if the storage of factors is affordable, then use the full AL preconditioner as the most robust and free of parameter tuning; otherwise switch to modified AL and adjust $\gamma$ to achieve iteration numbers somewhat higher but comparable to the full AL case.}

\section{Conclusions}
\rev{We conclude that recycling  AL preconditioner for the Oseen problem over several time steps is a highly effective strategy to reduce linear algebra costs when solving the time-dependent Navier--Stokes equations. In contrast to time-lagging of inertia terms, this approach does not compromise the numerical stability of the scheme.
For many 2D flows it is feasible to compute exact factorizations of velocity sub-blocks. Within the developed  framework  this ensures the solver robustness with respect to the Reynolds number. The performance of the approach was illustrated for a few examples of surface flows. The efficiency evident in numerical tests were backed by eigenvalue analysis, which extends some known results to the case of FE-level stabilization. To make the method even more efficient for smaller time steps, we introduced the pressure Schur complement preconditioner, which extends the Cahouet--Chabard preconditioner for the case of non-zero $(2,2)$-block $\vect C$. This extension can be useful in other settings with unfitted or pressure stabilized FEs.}

We expect that recycled modified AL preconditioner with a threshold ILU factorization for each sub-block can be an efficient strategy for large scale 3D flow problems. We plan to explore such possibility in the future. 

\section{Code availability}
\rev{The source code to run Kelvin--Helmholtz simulations from section~\ref{s_num} with building instructions and files defining input parameters was archived in~\cite{code}.}

\section{Acknowledgment}
This study does not have any conflicts to disclose.

\bibliographystyle{wileyNJD-VANCOUVER}
\bibliography{literatur}{}

\begin{thebibliography}{10}

\bibitem{benzi2006augmented}
Benzi M, and Olshanskii MA.
\newblock An augmented {L}agrangian-based approach to the {Oseen} problem.
\newblock SIAM Journal on Scientific Computing. 2006;{\bf 28}(6):2095--2113.

\bibitem{benzi2011modified}
Benzi M, Olshanskii MA, and Wang Z.
\newblock Modified augmented {L}agrangian preconditioners for the
  incompressible {N}avier--{S}tokes equations.
\newblock International Journal for Numerical Methods in Fluids. 2011;{\bf
  66}(4):486--508.

\bibitem{borm2012}
B{\"o}rm S, and Le~Borne S.
\newblock $\mathcal{H}$-{LU} factorization in preconditioners for augmented
  {L}agrangian and grad-div stabilized saddle point systems.
\newblock International journal for numerical methods in fluids. 2012;{\bf
  68}(1):83--98.

\bibitem{de2007two}
De~Niet A, and Wubs F.
\newblock Two preconditioners for saddle point problems in fluid flows.
\newblock International Journal for Numerical Methods in Fluids. 2007;{\bf
  54}(4):355--377.

\bibitem{farrell2020augmented}
Farrell PE, and Gazca-Orozco PA.
\newblock An Augmented Lagrangian Preconditioner for Implicitly Constituted
  Non-Newtonian Incompressible Flow.
\newblock SIAM Journal on Scientific Computing. 2020;{\bf 42}(6):B1329--B1349.

\bibitem{he2011augmented}
He X, Neytcheva M, and Capizzano SS.
\newblock On an augmented {L}agrangian-based preconditioning of {Oseen} type
  problems.
\newblock BIT Numerical Mathematics. 2011;{\bf 51}(4):865--888.

\bibitem{he2018combining}
He X, Vuik C, and Klaij CM.
\newblock Combining the augmented {L}agrangian preconditioner with the simple
  {Schur} complement approximation.
\newblock SIAM Journal on Scientific Computing. 2018;{\bf 40}(3):A1362--A1385.

\bibitem{heister2013efficient}
Heister T, and Rapin G.
\newblock Efficient augmented {{L}agrangian-type} preconditioning for the
  {Oseen} problem using {Grad-Div} stabilization.
\newblock International Journal for Numerical Methods in Fluids. 2013;{\bf
  71}(1):118--134.

\bibitem{moulin2019augmented}
Moulin J, Jolivet P, and Marquet O.
\newblock Augmented {L}agrangian preconditioner for large-scale hydrodynamic
  stability analysis.
\newblock Computer Methods in Applied Mechanics and Engineering. 2019;{\bf
  351}:718--743.

\bibitem{olshanskii2008augmented}
Olshanskii MA, and Benzi M.
\newblock An augmented {L}agrangian approach to linearized problems in
  hydrodynamic stability.
\newblock SIAM Journal on Scientific Computing. 2008;{\bf 30}(3):1459--1473.

\bibitem{ur2008comparison}
ur~Rehman M, Vuik C, and Segal G.
\newblock A comparison of preconditioners for incompressible {N}avier--{S}tokes
  solvers.
\newblock International Journal for Numerical Methods in Fluids. 2008;{\bf
  57}(12):1731--1751.

\bibitem{farrell2020reynolds}
Farrell PE, Mitchell L, Scott LR, and Wechsung F.
\newblock A {R}eynolds-robust preconditioner for the {R}eynolds-robust
  {S}cott--{V}ogelius discretization of the stationary incompressible
  {N}avier--{S}tokes equations.
\newblock arXiv preprint arXiv:200409398. 2020;.

\bibitem{farrell2019augmented}
Farrell PE, Mitchell L, and Wechsung F.
\newblock An {A}ugmented {L}agrangian Preconditioner for the 3{D} Stationary
  Incompressible {N}avier--{S}tokes Equations at High {R}eynolds Number.
\newblock SIAM Journal on Scientific Computing. 2019;{\bf 41}(5):A3073--A3096.

\bibitem{dahl1992ilu}
Dahl O, and Wille S.
\newblock An {ILU} preconditioner with coupled node fill-in for iterative
  solution of the mixed finite element formulation of the {2D} and {3D}
  {N}avier--{S}tokes equations.
\newblock International Journal for Numerical Methods in Fluids. 1992;{\bf
  15}(5):525--544.

\bibitem{konshin2017lu}
Konshin I, Olshanskii M, and Vassilevski Y.
\newblock {LU} factorizations and {ILU} preconditioning for stabilized
  discretizations of incompressible {N}avier--{S}tokes equations.
\newblock Numerical Linear Algebra with Applications. 2017;{\bf 24}(3):e2085.

\bibitem{konshin2015ilu}
Konshin IN, Olshanskii MA, and Vassilevski YV.
\newblock {ILU} Preconditioners for Nonsymmetric Saddle-Point Matrices with
  Application to the Incompressible {N}avier--{S}tokes Equations.
\newblock SIAM Journal on Scientific Computing. 2015;{\bf 37}(5):A2171--A2197.

\bibitem{segal2010preconditioners}
Segal A, Ur~Rehman M, and Vuik C.
\newblock Preconditioners for incompressible {N}avier--{S}tokes solvers.
\newblock Numerical Mathematics: Theory, Methods and Applications. 2010;{\bf
  3}(3):245--275.

\bibitem{bonito2020divergence}
Bonito A, Demlow A, and Licht M.
\newblock A divergence-conforming finite element method for the surface
  {Stokes} equation.
\newblock SIAM Journal on Numerical Analysis. 2020;{\bf 58}(5):2764--2798.

\bibitem{brandner2020finite}
Brandner P, and Reusken A.
\newblock Finite element error analysis of surface {Stokes} equations in stream
  function formulation.
\newblock ESAIM: Mathematical Modelling and Numerical Analysis. 2020;{\bf
  54}(6):2069--2097.

\bibitem{fries2018higher}
Fries TP.
\newblock Higher-order surface {FEM} for incompressible {Navier--Stokes} flows
  on manifolds.
\newblock International Journal for Numerical Methods in Fluids. 2018;{\bf
  88}(2):55--78.

\bibitem{gross2019meshfree}
Gross BJ, Trask N, Kuberry P, and Atzberger PJ.
\newblock Meshfree methods on manifolds for hydrodynamic flows on curved
  surfaces: a generalized moving least-squares {(GMLS)} approach.
\newblock Journal of Computational Physics. 2020;{\bf 409}:109340.

\bibitem{Jankuhn1}
Jankuhn T, Olshanskii MA, and Reusken A.
\newblock Incompressible Fluid Problems on Embedded Surfaces: Modeling and
  Variational Formulations.
\newblock Interfaces and Free Boundaries. 2018;{\bf 20}:353--377.

\bibitem{jankuhn2019higher}
Jankuhn T, and Reusken A.
\newblock Higher order trace finite element methods for the surface {Stokes}
  equation.
\newblock Preprint arXiv:190908327. 2019;.

\bibitem{lederer2020divergence}
Lederer PL, Lehrenfeld C, and Sch{\"o}berl J.
\newblock Divergence-free tangential finite element methods for incompressible
  flows on surfaces.
\newblock International Journal for Numerical Methods in Engineering. 2020;{\bf
  121}(11):2503--2533.

\bibitem{nitschke2019hydrodynamic}
Nitschke I, Reuther S, and Voigt A.
\newblock Hydrodynamic interactions in polar liquid crystals on evolving
  surfaces.
\newblock Physical Review Fluids. 2019;{\bf 4}(4):044002.

\bibitem{nitschke2012finite}
Nitschke I, Voigt A, and Wensch J.
\newblock A finite element approach to incompressible two-phase flow on
  manifolds.
\newblock Journal of Fluid Mechanics. 2012;{\bf 708}:418--438.

\bibitem{olshanskii2018finite}
Olshanskii MA, Quaini A, Reusken A, and Yushutin V.
\newblock A finite element method for the surface {Stokes} problem.
\newblock SIAM Journal on Scientific Computing. 2018;{\bf 40}(4):A2492--A2518.

\bibitem{olshanskiiinf}
Olshanskii MA, Reusken A, and Zhiliakov A.
\newblock Inf-sup stability of the trace {P2-P1} {T}aylor--{H}ood elements for
  surface {PDEs}.
\newblock Mathematics of Computation. 2020;.

\bibitem{olshanskii2019penalty}
Olshanskii MA, and Yushutin V.
\newblock A Penalty Finite Element Method for a Fluid System Posed on Embedded
  Surface.
\newblock Journal of Mathematical Fluid Mechanics. 2019;{\bf 21}(1):14.

\bibitem{reuther2018solving}
Reuther S, and Voigt A.
\newblock Solving the incompressible surface {Navier--Stokes} equation by
  surface finite elements.
\newblock Physics of Fluids. 2018;{\bf 30}(1):012107.

\bibitem{GurtinMurdoch75}
Gurtin ME, and Murdoch AI.
\newblock A continuum theory of elastic material surfaces.
\newblock Archive for Rational Mechanics and Analysis. 1975;{\bf
  57}(4):291--323.

\bibitem{rangamani2013interaction}
Rangamani P, Agrawal A, Mandadapu KK, Oster G, and Steigmann DJ.
\newblock Interaction between surface shape and intra-surface viscous flow on
  lipid membranes.
\newblock Biomechanics and modeling in mechanobiology. 2013;p. 1--13.

\bibitem{torres2019modelling}
Torres-S{\'a}nchez A, Mill{\'a}n D, and Arroyo M.
\newblock Modelling fluid deformable surfaces with an emphasis on biological
  interfaces.
\newblock Journal of Fluid Mechanics. 2019;{\bf 872}:218--271.

\bibitem{ORG09}
Olshanskii MA, Reusken A, and Grande J.
\newblock A Finite Element method for elliptic equations on surfaces.
\newblock SIAM J Numer Anal. 2009;{\bf 47}:3339--3358.

\bibitem{olshanskii2016trace}
Olshanskii MA, and Reusken A.
\newblock Trace Finite Element Methods for {PDEs} on Surfaces.
\newblock In: Bordas SPA, Burman E, Larson MG, and Olshanskii MA, editors.
  Geometrically Unfitted Finite Element Methods and Applications. Cham:
  Springer International Publishing; 2017. p. 211--258.

\bibitem{olshanskii2002low}
Olshanskii MA.
\newblock A low order {G}alerkin finite element method for the
  {N}avier--{S}tokes equations of steady incompressible flow: a stabilization
  issue and iterative methods.
\newblock Computer Methods in Applied Mechanics and Engineering. 2002;{\bf
  191}(47-48):5515--5536.

\bibitem{olshanskii2004grad}
Olshanskii M, and Reusken A.
\newblock Grad-div stablilization for {Stokes} equations.
\newblock Mathematics of Computation. 2004;{\bf 73}(248):1699--1718.

\bibitem{benzi2005numerical}
Benzi M, Golub GH, and Liesen J.
\newblock Numerical solution of saddle point problems.
\newblock Acta numerica. 2005;{\bf 14}:1.

\bibitem{ElmanBook}
Elman H, Silvester D, and Wathen A.
\newblock Finite Elements and Fast Iterative Solvers.
\newblock Oxford: Oxford University Press; 2005.

\bibitem{Temam84}
Temam R.
\newblock Navier--{Stokes} equations, theory and numerical analysis.
\newblock Amsterdam: North-Holland; 3rd ed., 1984.

\bibitem{GR}
Girault V, and Raviart PA.
\newblock Finite Element Methods for {Navier--Stokes} Equations.
\newblock Berlin: Springer; 1986.

\bibitem{burman2015stabilized}
Burman E, Hansbo P, and Larson MG.
\newblock A stabilized cut finite element method for partial differential
  equations on surfaces: The {L}aplace--{B}eltrami operator.
\newblock Computer Methods in Applied Mechanics and Engineering. 2015;{\bf
  285}:188--207.

\bibitem{burman2017cut}
Burman E, Hansbo P, Larson MG, and Massing A.
\newblock A cut discontinuous {G}alerkin method for the {L}aplace--{B}eltrami
  operator.
\newblock IMA Journal of Numerical Analysis. 2017;{\bf 37}(1):138--169.

\bibitem{burman2018cut}
Burman E, Hansbo P, Larson MG, and Massing A.
\newblock Cut finite element methods for partial differential equations on
  embedded manifolds of arbitrary codimensions.
\newblock ESAIM: Mathematical Modelling and Numerical Analysis. 2018;{\bf
  52}(6):2247--2282.

\bibitem{grande2017higher}
Grande J, Lehrenfeld C, and Reusken A.
\newblock Analysis of a High-Order Trace Finite Element Method for {PDEs} on
  Level Set Surfaces.
\newblock SIAM Journal on Numerical Analysis. 2018;{\bf 56}(1):228--255.
\newblock Available from: \url{https://doi.org/10.1137/16M1102203}.

\bibitem{olshanskii2009grad}
Olshanskii M, Lube G, Heister T, and L{\"o}we J.
\newblock Grad--div stabilization and subgrid pressure models for the
  incompressible {N}avier--{S}tokes equations.
\newblock Computer Methods in Applied Mechanics and Engineering. 2009;{\bf
  198}(49-52):3975--3988.

\bibitem{john2010numerical}
John V, and Kindl A.
\newblock Numerical studies of finite element variational multiscale methods
  for turbulent flow simulations.
\newblock Computer Methods in Applied Mechanics and Engineering. 2010;{\bf
  199}(13-16):841--852.

\bibitem{jankuhn2020error}
Jankuhn T, Olshanskii MA, Reusken A, and Zhiliakov A.
\newblock Error analysis of higher order trace finite element methods for the
  surface {S}tokes equations.
\newblock Journal of Numerical Mathematics. 2020;.

\bibitem{benzi2011analysis}
Benzi M, and Wang Z.
\newblock Analysis of augmented Lagrangian-based preconditioners for the steady
  incompressible Navier--Stokes equations.
\newblock SIAM Journal on Scientific Computing. 2011;{\bf 33}(5):2761--2784.

\bibitem{cahouet1988some}
Cahouet J, and Chabard JP.
\newblock Some fast 3{D} finite element solvers for the generalized {S}tokes
  problem.
\newblock International Journal for Numerical Methods in Fluids. 1988;{\bf
  8}(8):869--895.

\bibitem{schroeder2019reference}
Schroeder PW, John V, Lederer PL, Lehrenfeld C, Lube G, and Sch{\"o}berl J.
\newblock On reference solutions and the sensitivity of the {2D
  Kelvin--Helmholtz} instability problem.
\newblock Computers \& Mathematics with Applications. 2019;{\bf
  77}(4):1010--1028.

\bibitem{DROPS}
{DROPS package};.
\newblock \verb|http://www.igpm.rwth-aachen.de/DROPS/|.

\bibitem{trilinos-website}
{T}rilinos~{P}roject {T}eam T. The {T}rilinos {P}roject {W}ebsite;.

\bibitem{lehrenfeld2016high}
Lehrenfeld C.
\newblock High order unfitted finite element methods on level set domains using
  isoparametric mappings.
\newblock Computer Methods in Applied Mechanics and Engineering. 2016;{\bf
  300}:716--733.

\bibitem{code}
DROPS: surface Navier-Stokes solver.
\newblock
  \url{https://github.com/56th/drops/archive/refs/heads/surfaceNSE_06/03/2021.zip};
  2021.

\bibitem{bendixson1902racines}
Bendixson I.
\newblock Sur les racines d'une {\'e}quation fondamentale.
\newblock Acta Mathematica. 1902;{\bf 25}(1):359--365.

\end{thebibliography}


\appendix

\section{Proof of (12)}
Since the grad-div stabilization does not deliver the algebraic structure of the  augmented Lagrangian  as in \cite{benzi2006augmented,benzi2011modified},  we cannot make use of the Sherman–Morrison–Woodbury formula or similar representations of the pressure Schur complement of the augmented system.  Therefore we base our proof of \eqref{bound} on a different argument. Let $\mathbb{R}^m$ and $\mathbb{R}^n$ be the coefficient spaces for the pressure and velocity finite element functions, respectively. For $\mathrm p\in \mathbb{R}^m$ the corresponding finite element function is $p_h\in Q_h$, similar we have $\bu_h\in \bV_h$ for $\mathrm u\in \mathbb{R}^n$, etc. Further $\la\cdot,\cdot\ra$ denotes the Euclidian inner product and $\|\cdot\| \coloneqq \la\cdot,\cdot\ra^{\frac12}$.

The low bound for the real parts of the eigenvalues from \eqref{eig_pr} is given by the Bendixson theorem~\cite{bendixson1902racines}:
\[
\inf_{\mathrm p\in \mathbb{R}^m}\frac{\la\frac12(\bS+\bS^T)\,\mathrm p, \mathrm p\ra}{\la \bQ\, \mathrm p, \mathrm p\ra}\le \Re(\lambda), \quad\text{with } \bQ=(\nu + \gamma)^{-1}\,\vect M_p + \vect C.
\]
Noting that $\la\frac12(\bS+\bS^T)\, \mathrm p, \mathrm p\ra = \la\bS\, \mathrm p, \mathrm p\ra = \la \bB\, \mathrm u, \mathrm p\ra + \la \bC\,\mathrm p,\mathrm p\ra$, $\bA\,\mathrm u = \bB^T \mathrm p$, we re-write the quantities on the left-hand side in the finite element notation:
\begin{equation}\label{aux671}
\begin{split}
\la\bS\,\mathrm p, \mathrm p\ra &= b(\bu_h,p_h)+s(p_h,p_h)\\
b(\bv_h,p_h)& =a(\bu_h,\bv_h)+\gamma(\bu_h,\bv_h),\quad \forall~\bv_h\in\bV_h.
\end{split}
\end{equation}
Letting $\bv_h=\bu_h$  and using the skew-symmetric form of the advection term: 
\begin{equation}\label{aux691}
\frac12\int_\Gamma\left(\bv_h^T(\nablaG\bu_h)\bw- \bu_h^T(\nablaG\bv_h)\bw\right)\,ds,
\end{equation} 
we get
\begin{equation}\label{aux690}
\la\bS\,\mathrm p,\mathrm p\ra= a(\bu_h,\bu_h)+\gamma \|\divG\bu_h\|_{L^2(\Gamma)}^2+s(p_h,p_h)
\end{equation}
with $\bu_h$ solving the second equation in \eqref{aux671} for the given $p_h$, i.e. $\mathrm u = \bA^{-1}\bB^T \mathrm p$.

Let $\|\bv_h\|_\ast^2 \coloneqq 2\int_\Gamma  |E_s( \bv_h)|^2\, ds +\tau \int_{\Gamma}|\bn\cdot\bv_h|^2 \,ds + \rho_u \int_{\OGamma} |(\bn\cdot\nabla) \bv_h|^2  \, dx$.
The inf-sup condition  for trace $\vect P_2$--$P_1$ elements proved in~\cite{olshanskiiinf} reads:
\begin{equation}\label{LBBb}
 c_0^2\,\|p_h\|_{L^2(\Gamma)}^2
  \le \sup_{\bv_h\in\bV_h}\frac{b(\bv_h,p_h)^2}{\|\bv_h\|_\ast^2} + s(p_h,p_h),
\end{equation}
with  $c_0>0$ independent of $h$ and position of $\Gamma$ in the background mesh.  The condition \eqref{LBBb} can be rewritten as follows: There exists $\bv_h\in\bV_h$  such that
\[
c_0^2\,\|p_h\|_{L^2(\Gamma)}^2 \le b(\bv_h,p_h)^2 + s(p_h,p_h),\quad\text{and}~\|\bv_h\|_\ast=1.
\]
We take this $\bv_h$ as a test function in \eqref{aux671} and apply the Cauchy--Schwarz and Poincar\'{e}  inequalities to arrive at
\[
\begin{split}
c_0^2\,\|p_h\|_{L^2(\Gamma)}^2
	&\le [a(\bu_h,\bv_h)+\gamma(\bu_h,\bv_h)]^2 + s(p_h,p_h)\\
	&\le [a(\bu_h,\bu_h)^{\frac12}a(\bv_h,\bv_h)^{\frac12}+\gamma\|\divG\bu_h\|_{L^2(\Gamma)}\|\divG\bv_h\|_{L^2(\Gamma)} + \\
	&\quad\:\,\|\bw\|_{L^\infty(\Gamma)}\|\nabla\bu_h\|_{L^2(\Gamma)}\|\nabla\bv_h\|_{L^2(\Gamma)}]^2 + s(p_h,p_h)\\
	&\le [a(\bu_h,\bu_h)^{\frac12}+\gamma\|\divG\bu_h\|_{L^2(\Gamma)}+C\,\|E_s(\bu_h)\|_{L^2(\Gamma)}]^2 + s(p_h,p_h)
\end{split}
\]
where in the last inequality we used
\[ a(\bv_h,\bv_h)\le \|\bv_h\|^2_\ast=1,\quad \|\divG\bv_h\|_{L^2(\Gamma)}\le\|E_s(\bv_h)\|_{L^2(\Gamma)}\le \|\bv_h\|_\ast=1,\quad \|\bw\|_{L^\infty(\Gamma)}=1,\]
and the surface Korn inequality.  Since $\|E_s(\bu_h)\|_{L^2(\Gamma)}^2\le \nu^{-1}a(\bu_h,\bu_h)$ we obtain thanks to \eqref{aux690} the estimate
\[
\begin{split}
c_0^2\,\|p_h\|_{L^2(\Gamma)}^2 &\le 3[(1+C\,\nu^{-1})a(\bu_h,\bu_h) + \gamma^2 \|\divG\bu_h\|_{L^2(\Gamma)}^2]+s(p_h,p_h)\\
&\le
3(1+C\,\nu^{-1}+\gamma)\la\bS\,\mathrm p, \mathrm p\ra.
\end{split}
\]
Since $\|p_h\|_{L^2(\Gamma)}^2=\la M_p\,\mathrm p, \mathrm p\ra$, we get 
\[
\frac{c_0^2(\nu+\gamma)}{3(1+C\,\nu^{-1}+\gamma)} \la(\nu+\gamma)^{-1} \vect M_p\,\mathrm p, \mathrm p\ra\le
\la\bS\,\mathrm p, \mathrm p\ra.
\]
Finally, using $\nu\le 1$ and $\la\bC\,\mathrm p, \mathrm p\ra\le \la\bS\,\mathrm p, \mathrm p\ra$ the above estimate yields
\[
\frac{c\,(\nu+\gamma)}{\nu^{-1}+\gamma}\left( \la(\nu+\gamma)^{-1} \vect M_p\,\mathrm p, \mathrm p\ra+\la\bC\,\mathrm p, \mathrm p\ra\right) \le
\la\bS\,\mathrm p, \mathrm p\ra+\la\bC\,\mathrm p, \mathrm p\ra.
\]
with some $c$ independent of $h$ and position of $\Gamma$ in the background mesh.
To show the bound on $|\lambda|$, we estimate
\[
|\lambda|\le \|\bQ^{-\frac12}\bS\bQ^{-\frac12}\|=\sup_{\mathrm p, \mathrm q\in \mathbb{R}^m}\frac{\la\bQ^{-\frac12}\bS\bQ^{-\frac12}\,\mathrm p, \mathrm q\ra}{\| \mathrm p\|\|\mathrm q\|}
=\sup_{\mathrm p, \mathrm q\in \mathbb{R}^m}\frac{\la\bS\,\mathrm p, \mathrm q\ra}{\|\bQ^{\frac12}\,\mathrm p\|\|\bQ^{\frac12}\,\mathrm q\|}.
\]
In finite element notations, we rewrite
\begin{equation}\label{aux728}
\begin{split}
\la\bS\,\mathrm p, \mathrm q\ra &= b(\bu_h,q_h)+s(p_h,q_h)\\
b(\bv_h,p_h)  &=a(\bu_h,\bv_h)+\gamma(\bu_h,\bv_h),\quad \forall~\bv_h\in\bV_h.
\end{split}
\end{equation}
The Cauchy--Schwarz inequality yields
\[
\la\bS\,\mathrm p, \mathrm q\ra \le \|\divG\bu_h\|_{L^2(\Gamma)}\|q_h\|_{L^2(\Gamma)}+s(p_h,p_h)^{\frac12}s(q_h,q_h)^{\frac12}
\]
and
\[
a(\bu_h,\bu_h)+\gamma\|\divG\bu_h\|_{L^2(\Gamma)}^2=b(\bu_h,p_h)\le\|\divG\bu_h\|_{L^2(\Gamma)}\|p_h\|_{L^2(\Gamma)}.
\]
From the second equation we get $\gamma\|\divG\bu_h\|_{L^2(\Gamma)}\le\|p_h\|_{L^2(\Gamma)}$, and substituting this into the first equation we get
\[
\la\bS\,\mathrm p, \mathrm q\ra \le\frac{\nu+\gamma}{\gamma} \left\|\frac1{\sqrt{\nu+\gamma}}p_h\right\|_{L^2(\Gamma)}\left\|\frac1{\sqrt{\nu+\gamma}}q_h\right\|_{L^2(\Gamma)}+s(p_h,p_h)^{\frac12}s(q_h,q_h)^{\frac12}
\le
\frac{\nu+\gamma}{\gamma} \|\bQ^{\frac12}\,\mathrm p\|\|\bQ^{\frac12}\,\mathrm q\|.
\]
This proves the desired bounds.

\end{document}